\documentclass[a4paper,11pt]{article}
\usepackage{amsfonts}
\usepackage{amssymb}
\usepackage{amsthm}
\usepackage{pstricks}
\usepackage{amsmath}
\usepackage{mathrsfs}

\def\XXint#1#2#3{{\setbox0=\hbox{$#1{#2#3}{\int}$}
     \vcenter{\hbox{$#2#3$}}\kern-.5\wd0}}

\theoremstyle{definition} \theoremstyle{remark}
\newtheorem{thm}{Theorem}
\newtheorem{exm}{Example}
\newtheorem{lmm}{Lemma}
\newtheorem{dfn}{Definition}
\newtheorem{rmk}{Remark}
\newtheorem{clr}{Corollary}

\DeclareMathOperator{\dist}{dist}
\DeclareMathOperator{\diag}{diag}

\begin{document}

\title{\textbf{Fractional Calculus in $\mathbb{C}$ and linear ODE's with
analytic coefficients. Part 1}}
\author{\textbf{V.P. Gurarii}}
\maketitle

\begin{abstract}
The theory of fractional calculus in the complex plane was not built with a
specific application in mind. The main obstacle to application was the
difficulty with obtaining analytic continuations of fractional derivatives
and integrals. It is known that if a function is analytic in a simply
connected domain, then its fractional derivatives and integrals are also
analytic in this region. In multi-connected domains, in general, this
property does not hold. However, for the set of hypergeometric functions $%
_{p+1}F_{p},p\geq 2,$ this property is preserved. We propose a new version
of fractional calculus, which allows us to reveal the reason for this
preservation, and introduce broad classes of functions for which an
appropriate theory can be constructed. This allows us to find a meaningful
application of our calculus in the asymptotic theory of linear ODE's with
analytic coefficients.
\end{abstract}

\tableofcontents


\section{Contents of Part I and Part II}

\underline{\bf Part 1} \\ I. Fractional Calculus in $\mathbb{C}.$ \\ II. Linear ODE's with
analytic coefficients. \\ \underline{\bf Part 2}\\ III. Borel-type summation. \\ IV. Error bounds
and Stokes' phenomenon. \\ V. Connection coefficients.

\section{Fractional Calculus in $\mathbb{C}$}

\subsection{Introduction}

Fractional Calculus has existed for more than 300 years. The encyclopedic
monograph \cite{Samko} provides an insight into the progress made in this
area up to the early 1990s. In recent years the topic has been rediscovered
by scientists and engineers in an increasing number of fields in which
complex analysis is essential, see, for example, \cite{Or}--\cite{Das}. It
has been shown, see \cite{Osler1}, \cite{Osler2}, that if a function is
analytic in a simply connected domain then its fractional derivatives and
integrals are also analytic in this region. However, in the case of finitely
connected domains this property is, in general, not preserved. We have
constructed a version of fractional calculus in the complex plane that
enables us to understand more clearly the reason for the obstacle to
analytical continuation and to overcome the above obstacle for broad classes
of analytic functions suggested by linear ODE's with analytic coefficients.
We apply then our version of fractional calculus to the asymptotic analysis
of solutions of these ODE's. In the current paper we extend the results of
\cite{GG-1}, which are valid for ODE's with coefficients given by even
analytic functions, to the general case of the second order linear ODE's
with analytic coefficients.

In the next section, we introduce the two functional linear spaces in which
our version of fractional calculus is defined.

\subsection{Laplace-Borel dual spaces of analytic functions}

Let $H$ be the linear space that consists of all functions $F\left( t\right)
$ possessing the following properties:

\begin{itemize}
\item (i) For some $a>0$, $F\left( t\right) $ is analytic in the region
\begin{equation}
\mathfrak{D}(a)=\{t\in \mathbb{C}:\dist\{t,(0,+\infty )\}<a\}
\label{N-region}
\end{equation}%
with boundary
\begin{equation}
\gamma (a)=\{t\in \mathbb{C}:\dist\{t,(0,+\infty )=a\}.  \label{N-boundary}
\end{equation}

\item (ii) For some $R\geq 0,$ $F\left( t\right) $ has an exponential growth
of type $R$ at infinity of the region $\mathfrak{D}(a)$.
\end{itemize}

Given $F\left( t\right) \in H$ set
\begin{equation}
a\left( F\right) =\sup \{a:(i)\text{ is satisfied}\}.  \label{aF}
\end{equation}%
Likewise, we set%
\begin{equation}
R\left( F\right) =\inf \{R:(ii)\text{ is satisfied}\}.  \label{RF}
\end{equation}

We denote by $\mathcal{L}$ the Laplace transform operator of the form%
\begin{equation}
P\left( \zeta \right) =\mathcal{L}\left\{ F\left( t\right) \right\} =\zeta
\int_{0}^{+\infty }e^{-\zeta t}F\left( t\right) dt,  \label{LaplaceF}
\end{equation}%
where the integral (\ref{LaplaceF}) is absolutely convergent for any $\zeta
,\Re \{\zeta \}>R\geq 0.$

\begin{rmk}
In paper \cite{Gur}, where our version of fractional calculus was first
proposed, we denote by $\mathcal{L}$ the standard Laplace transform
operator, see equations (44) and (46) with $\alpha =0$. We note that keeping
in mind the application to differential equations with analytic
coefficients, it is a little more convenient to use (\ref{LaplaceF}) instead
of the standard definition. Note that all the definitions and the results
obtained here can be easily rewritten in terms of the standard Laplace
transform, and vice versa, using the change $P\rightarrow {\zeta P}$..
\end{rmk}

We set%
\begin{equation*}
S=\mathcal{L}H=\left\{ P:P\left( \zeta \right) =\mathcal{L}\left\{ F\left(
t\right) \right\} \right\} ,
\end{equation*}
where $\mathcal{L}$ is given by (\ref{LaplaceF}), $F\left( t\right) \in H,$
and the integral (\ref{LaplaceF}) is absolutely convergent for any $\zeta
,\Re \{\zeta \}>R\left( F\right) .$

In order to simplify our development, and to provide the space $H$ with a
countable normalized topology, we introduce subspaces $H(a,R)\subset H$ and $%
S(a,R)\subset S$.

\begin{itemize}
\item (i) If, in the definition of $H$, we fix the values of $a$ and $R$ we
obtain the definition of a subspace of $H$, which we denote by $H(a,R)$.

\item (ii) Now let $S(a,R)=\mathcal{L}H(a,R)$--a subspace of $S$.
\end{itemize}

Thus $S(a,R)$ and $H(a,R)$ are Laplace-Borel dual linear spaces and we have
\begin{equation*}
H=\cup _{a>0,R\geq 0}H(a,R),S=\cup _{a>0,R\geq 0}S(a,R).
\end{equation*}

It follows immediately that $\mathcal{L}:H\rightarrow S$ is a bijection.

Later we will present an intrinsic definition of $S$ and $S(a,R)$ in terms
of properties of their elements.

Since $F\left( t\right) \in H(a,R)$ is analytic in the circle $|t|<a$, its
Taylor series expansion is absolutely convergent inside this circle and we
have
\begin{equation}
F\left( t\right) =\sum_{k=0}^{\infty }f_{k}t^{k},f_{k}=\frac{F^{\left(
k\right) }\left( 0\right) }{k!}.  \label{Taylor-s}
\end{equation}

Suppose that the function $P\left( \zeta \right) =\mathcal{L}\left\{ F\left(
t\right) \right\} \in S(a,R).$ Then $P\left( \zeta \right) $ is analytic and
bounded in any half-plane $\Re \{\zeta \}>r>R=R\left( F\right) \geq 0$.
Furthermore, setting
\begin{equation}
p_{k}=f_{k}k!,k=0,1,\ldots ,  \label{Borel}
\end{equation}%
where $f_{k}$ is given by (\ref{Taylor-s}), we may approximate $P\left(
\zeta \right) $ by the partial sums $\sum_{k=0}^{n-1}p_{k}/\zeta ^{k}$ of
the power series $\sum_{k=0}^{\infty }p_{k}/\zeta ^{k}.$ We introduce the
remainder $P_{n}\left( \zeta \right) $ as
\begin{equation}
P_{n}\left( \zeta \right) =P\left( \zeta \right) -\sum_{k=0}^{n-1}\frac{p_{k}%
}{\zeta ^{k}},n=0,1,\ldots ,  \label{P-remainder}
\end{equation}%
where $P_{0}\left( \zeta \right) =P\left( \zeta \right) $. Clearly, $%
P_{n}\left( \zeta \right) \in S(a,R)$ and it follows from Watson's lemma
that for $\Re \left\{ \zeta \right\} >R$ the following asymptotic expansion
is valid
\begin{equation}
P_{n}\left( \zeta \right) \sim \sum_{k=n}^{\infty }\frac{p_{k}}{\zeta ^{k}}%
,\zeta \rightarrow \infty ,  \label{P-expansion}
\end{equation}%
which means that for $n=0,1,\ldots ,$ we have%
\begin{equation}
P_{n}\left( \zeta \right) =O\left( \frac{1}{\left\vert \zeta \right\vert ^{n}%
}\right) ,\zeta \rightarrow \infty .  \label{P-n- error}
\end{equation}%
In what follows we will work also with a function $F\left( t,\alpha \right)
\in H(a,R)$, generated by $F\left( t\right) $, which depends on a parameter $%
\alpha $ such that $F\left( t,0\right) =F\left( t\right) $. Then all the
formulas (\ref{Taylor-s})--(\ref{P-n- error}) are retained for $F\left(
t,\alpha \right) $ with the corresponding changes replacing $f_{k},P\left(
\zeta \right) ,p_{k},P_{n}\left( \zeta \right) $ by $f_{k}\left( \alpha
\right) ,P\left( \zeta ,\alpha \right) ,p_{k}\left( \alpha \right)
,P_{n}\left( \zeta ,\alpha \right),$, respectively. Our main aim is to find
a method of summation of the asymptotic series given by (\ref{P-expansion})
with $p_{k}=p_{k}\left( \alpha \right) $ that would enable us to represent
the results of summation of the series $\sum_{k=n}^{\infty }p_{k}\left(
\alpha \right) /\zeta ^{k}$ in the form of an integral transformation of $%
F\left( t,\alpha \right) .$ The required function $F\left( t,\alpha \right) $
and the kernel of this transformation have been presented, with the
appropriate change of notation, in Section 4 and Section 5 of \cite{Gur},
see equations (55), (100) and (105) where $s=\alpha +n$. Later we show how
to apply this technique to the asymptotic analysis of solutions of linear
ODE's with analytic coefficients.

It can be easily proved that functions $P\left( \zeta \right) \in S$ are
uniquely determined by the coefficients $p_{k}$ of their asymptotic
expansions. In other words, we claim the following statement.

\begin{thm}
{(Uniqueness Theorem)} If $P_{1}\left( \zeta \right) ,P_{2}\left( \zeta
\right) \in S$ and $p_{1,k}=p_{2,k}$ then $P_{1}\left( \zeta \right) \equiv
P_{2}\left( \zeta \right) $.
\end{thm}

Proof. Indeed, suppose $P_{1}\left( \zeta \right) =\mathcal{L}\left\{
F_{1}\left( t\right) \right\} $ and $P_{2}\left( \zeta \right) =\mathcal{L}%
\left\{ F_{2}\left( t\right) \right\} .$ Setting $a=\min (a(F_{1}),a(F_{2}))$
and $R=\max (R(F_{1}),R(F_{2})),$ it follows that $F_{1}\left( t\right) $
and $F_{2}\left( t\right) $ belong to $H(a,R)$. Thus using (\ref{Borel}), $%
p_{1,k}=p_{2,k}\Longrightarrow f_{1,k}=f_{2,k}\Longrightarrow F_{1}\left(
t\right) \equiv F_{2}\left( t\right) \Longrightarrow P_{1}\left( \zeta
\right) \equiv P_{2}\left( \zeta \right) .$

In the next section we introduce operators of fractional differentiation and
integration in the complex plane in the spaces $H(a,R)$ and $S(a,R)$ and
explain why we needed new definitions that are different from the classical
definitions.

\subsection{ Fractional derivatives and integrals in $\mathbb{C}$}

For complex $\alpha ,\Re \{\alpha \}>-1$, we introduce a family of operators
$\mathcal{L}_{\alpha }$ acting on $H$:%
\begin{equation}
\mathcal{L}_{\alpha }\left\{ F\left( t\right) \right\} =\zeta ^{\alpha }%
\mathcal{L}\left\{ t^{\alpha }F\left( t\right) \right\} ,  \label{L-alpha}
\end{equation}

For $\alpha =0$ we have $\mathcal{L}_{0}=\mathcal{L}$. The pair of operators
$\mathcal{L}$ and $\mathcal{L}_{\alpha }$ generates four linear operators:

\begin{align}
& \mathcal{D}_{\alpha }=\mathcal{\mathcal{L}}^{-1}\mathcal{L}_{\alpha },%
\text{ }\mathcal{I}_{\alpha }=\mathcal{\mathcal{L}}_{\alpha }^{-1}\mathcal{L}%
,  \label{inH} \\
& \mathcal{\hat{D}}_{\alpha }=\mathcal{L}_{\alpha }\mathcal{\mathcal{L}}%
^{-1},\,\mathcal{\hat{I}}_{\alpha }=\mathcal{L\mathcal{L}}_{\alpha }^{-1}.
\label{inS}
\end{align}

\begin{dfn}
The operators $\mathcal{D}_{\alpha }$ and $\mathcal{I}_{\alpha }$ are said
to be operators of fractional differentiation and integration of order $%
\alpha $\ acting on the space $H$. The second pair of operators $\mathcal{%
\hat{D}}_{\alpha }$ and $\mathcal{\hat{I}}_{\alpha }$ are the dual operators
of fractional differentiation and integration acting on the space $S$.
\end{dfn}

\begin{exm}
For any $a>0$ and a non-negative integer $k,$ $F\left( t\right) =t^{k}\in
H(a,0)$. Then for $\Re \{\alpha \}>-1$%
\begin{equation*}
\mathcal{L}_{\alpha }\left\{ t^{k}\right\} =\zeta ^{1+\alpha
}\int_{0}^{+\infty }e^{-\zeta t}t^{\alpha +k}dt=\frac{\Gamma \left( \alpha
+k+1\right) }{\zeta ^{k}},
\end{equation*}%
so that
\begin{equation}
\mathcal{D}_{\alpha }\left\{ t^{k}\right\} =\frac{\Gamma \left( \alpha
+k+1\right) }{k!}t^{k}.  \label{Dt^k}
\end{equation}%
Since $\mathcal{I}_{\alpha }=\mathcal{D}_{\alpha }^{-1},$%
\begin{equation}
\mathcal{I}_{\alpha }\left\{ t^{k}\right\} =\frac{k!}{\Gamma \left( \alpha
+k+1\right) }t^{k}.  \label{It^k}
\end{equation}
\end{exm}

Let us introduce the following notation
\begin{equation}
F\left( t,\alpha \right) =\mathcal{D}_{\alpha }\left\{ F\left( t\right)
\right\} ,F^{\ast }\left( t,\alpha \right) =\mathcal{I}_{\alpha }\left\{
F\left( t\right) \right\} .  \label{F-calculus}
\end{equation}

\begin{equation}
P\left( \zeta ,\alpha \right) =\mathcal{\hat{D}}_{\alpha }\left\{ P\left(
\zeta \right) \right\} ,P^{\ast }\left( \zeta ,\alpha \right) =\mathcal{\hat{%
I}}_{\alpha }\left\{ P\left( \zeta \right) \right\} .  \label{P-calculus}
\end{equation}%
It follows, using (\ref{Dt^k}) and (\ref{It^k}), that Taylor series of $%
F\left( t,\alpha \right) $ and $F^{\ast }\left( t,\alpha \right) $ {can be
written as}
\begin{equation}
F\left( t,\alpha \right) =\sum_{k=0}^{\infty }f_{k}\left( \alpha \right)
t^{k}\text{ and }F^{\ast }\left( t,\alpha \right) =\sum_{k=0}^{\infty
}f_{k}^{\ast }\left( \alpha \right) t^{k},  \label{T-ser}
\end{equation}%
{where}%
\begin{equation}
f_{k}\left( \alpha \right) =\frac{f_{k}\Gamma \left( \alpha +k+1\right) }{k!}%
\text{ and }f_{k}^{\ast }\left( \alpha \right) =\frac{f_{k}k!}{\Gamma \left(
\alpha +k+1\right) },  \label{T-coef}
\end{equation}%
respectively. Clearly, the expansions in (\ref{T-ser}) are absolutely
convergent {in the circle }$\left\vert t\right\vert <a$.

Using integral representations for $F\left( t,\alpha \right) $ and $F^{\ast
}\left( t,\alpha \right) $ in terms of $F(t)\in H(a,R)$, which will be
demonstrated later, we can claim the following results, valid for all $a>0,
R\geq 0$ and every $\alpha ,\Re \{\alpha \}>-1:$

\begin{thm}
\begin{itemize}
\item (i) $F\left( t,\alpha \right), F^{\ast}\left( t,\alpha \right)\in
H(a,R),$

\item (ii)
\begin{equation}
\mathcal{D}_{\alpha }H(a,R)=H(a,R),\mathcal{I}_{\alpha }H(a,R)=H(a,R),
\label{CinH}
\end{equation}%
\begin{equation}
\mathcal{\hat{D}}_{\alpha }S(a,R)=S(a,R),\mathcal{\hat{I}}_{\alpha
}S(a,R)=S(a,R),  \label{CinS}
\end{equation}%
\begin{equation}
\mathcal{D}_{\alpha }H=H,\mathcal{I}_{\alpha }H=H,  \label{CH}
\end{equation}%
\begin{equation}
\mathcal{\hat{D}}_{\alpha }S=S,\mathcal{\hat{I}}_{\alpha }S=S,  \label{SH}
\end{equation}%
{and every operator in (\ref{CinH})--(\ref{SH}) is bijective.}
\end{itemize}
\end{thm}

\begin{rmk}
\begin{itemize}
\item The function $F\left( t,\alpha \right) =\mathcal{D}_{\alpha }\left\{
F\left( t\right) \right\} $ was first introduced in Section 4, Eq.(55) of
\cite{Gur}. The main result of that Section is the derivation of the
integral representation for $F\left( t,\alpha \right) $ in terms of $F\left(
t\right) $. Further development of the technique from \cite{Gur} enables us
derive the integral representation for $F^{*}\left( t,\alpha \right) $ in
terms of $F\left(t\right),$ to prove the validity of (\ref{CinH}), and then
to study the problem of the analytic continuations of $\mathcal{D}_{\alpha
}\left\{ F\left( t\right) \right\} $ and $\mathcal{I}_{\alpha }\left\{
F\left( t\right) \right\} $ in the $t$ and $\alpha $-planes.

\item As far as we know, our definitions of fractional derivatives and
integrals given by (\ref{inH}) and (\ref{inS}) are new and can be extended
to the cases when the Laplace operators are replaced by the
Laplace-Stieltjes operators. This version is different from the classical
one. Of the five common standard requirements for different variations of
the classical fractional calculus given in \cite{Ross-1}, our version
retains only three. On the other hand, there is a link between our
definitions and the classical Liouville-Riemann definitions. There is also a
similarity between our version and the theory of Fractional Differ-Integrals
by Grunwald-Letnikov.

\item We note the following fact, which gives some support to the legitimacy
of our version of fractional calculus. For the particular case when $F\left(
t\right) \in H$ and $F\left( t\right) /t$ is absolutely integrable on the
contour given by (\ref{N-boundary}) the formulas giving the integral
representations for $\mathcal{D}_{\alpha }\left\{ F\left( t\right) \right\} $
and $\mathcal{I}_{\alpha }\left\{ F\left( t\right) \right\} $ are much
simpler than for the general case. They are given by (\ref{alpha-dual-D})
and (\ref{alpha-dual-I}) in Section 2.5. Similar formulas for a version of
fractional derivatives and integrals have appeared earlier in the monograph
\cite{Samko}, in the section Fractional calculus in the complex plane, the
only difference being that our path of integration given by (\ref{N-boundary}%
) replaces a circle centered at the origin of radius $a$ used there. We note
that our definitions were at that time unknown. We add, in conclusion, that
the construction of this version of fractional calculus is not an end in
itself for us. This theory proved to be necessary for the implementation of
our approach to the asymptotic analysis of linear ODE's with analytic
coefficients.
\end{itemize}
\end{rmk}

Let $F\left( t\right) \in H(a,R)$ and
\begin{equation}
P\left( \zeta ,\alpha \right) =\mathcal{L}_{\alpha }\left\{ F\left( t\right)
\right\} =\zeta ^{1+\alpha }\int_{0}^{+\infty }e^{-\zeta t}t^{\alpha
}F\left( t\right) dt,  \label{LM-tr}
\end{equation}
where $\mathcal{L}_{\alpha }$ is given by (\ref{L-alpha}). Then result (\ref%
{LM-tr}) yields the following statement:

\begin{thm}
\begin{equation}
P\left( \zeta ,\alpha \right) =\zeta \int_{0}^{+\infty }e^{-\zeta t}\mathcal{%
D}_{\alpha }\left\{ F\left( t\right) \right\} dt,  \label{LM-alt}
\end{equation}%
and for $\Re \{\zeta \}>R$ the integral (\ref{LM-alt}) is absolutely
convergent.
\end{thm}

Proof. Using relations (\ref{T-ser}) and (\ref{T-coef}), it can be verified
that the power series expansions in $1/\zeta $ of both the integrals in (\ref%
{LM-tr}) and (\ref{LM-alt}) coincide. The result (\ref{LM-tr}) then follows
from Theorem 1.

\begin{rmk}
We note that the Laplace-Mellin transforms of the form (\ref{LM-tr}) have
been known for a long time. However, the fact that there is an alternative
representation in the form of the Laplace transform (\ref{LM-alt}) was not
known except for the cases when $F\left( t\right) $ is the hypergeometric
function. See, for example, the integral representations of confluent
hypergeometric functions in \cite{Abramowitz}, \cite{Erdelyi} and \cite%
{Buchholz} in the form of the Laplace transform and Laplace-Mellin transform
of hypergeometric functions.
\end{rmk}

In the future, for applications to the asymptotic theory of differential
equations, we require the following result which follows from Theorem 1, and
gives, in fact, alternative definitions of fractional derivatives and
integrals in the space $H(a,R)$.

\begin{thm}
Let $F\left( t\right) \in H(a,R)$. Then
\begin{align}
&\zeta \int_{0}^{\infty }e^{-\zeta t}t^{\alpha }F\left( t\right) dt=\zeta
^{1-\alpha }\int_{0}^{\infty }e^{-\zeta t}\mathcal{D}_{\alpha }\left\{
F\left( t\right) \right\} dt,  \label{alt-D} \\
&\zeta \int_{0}^{\infty }e^{-\zeta t}F\left( t\right) dt=\zeta ^{1+\alpha
}\int_{0}^{\infty }e^{-\zeta t}t^{\alpha }\mathcal{I}_{\alpha }\left\{
F\left( t\right) \right\} dt,  \label{alt-I}
\end{align}%
and all integrals are absolutely convergent for $\Re \{\zeta \}>R$.
\end{thm}

In what follows we illustrate what has been said in this section by applying
the above formula to the case when $F\left( t\right) $ is a generalized
hypergeometric function. This case forms an important touchstone for our
work.

\subsection{Hypergeometric Functions}

Given $p=1,2,\ldots ,$ we consider the hypergeometric series of the form
\begin{equation}
\left( _{p+1}F_{p}\left( a_{1},\ldots ,a_{p+1};b_{1},\ldots ,b_{p};-t\right)
\right) =\sum_{k=0}^{\infty }\left( -1\right) ^{k}\frac{\left( a_{1}\right)
_{k}\ldots \left( a_{p+1}\right) _{k}}{\left( b_{1}\right) _{k}\ldots \left(
b_{p}\right) _{k}}\frac{t^{k}}{k!},  \label{H-p}
\end{equation}%
\ It is known, see \cite{Erdelyi}, that this series represents a
hypergeometric function $\left( _{p+1}F_{p}\right) \left( -t\right) $ which
is analytic in the unit circle of the $t$-plane and admits an analytical
continuation to the $t$-plane cut along the interval $\left( -\infty
,-1\right) $. Moreover, this function admits further analytic continuation
to the extended $t$-plane except for the points $0,-1,\infty $ which are
regular singularities. It follows that for $p=0,1,\ldots ,$ the function $%
\left( _{p+1}F_{p}\right) \left( -t\right) \ $belongs to $H\left( a,R\right)
$, where $a=1,R=0$, and\ using our definitions given by (\ref{inH}) and (\ref%
{inS}) and relations (\ref{T-ser}) and (\ref{T-coef}), we have%
\begin{align}
&\mathcal{D}_{\alpha }\left\{ \left( _{p+1}F_{p}\right) \left( -t\right)
\right\}  \notag \\
=&\Gamma \left( \alpha +1\right) \left( _{p+2}F_{p+1}\left( a_{1},\ldots
,a_{p+1},\alpha +1;b_{1},\ldots ,b_{p},1;-t\right) \right),  \label{D-H-p} \\
&\mathcal{I}_{\alpha }\left\{ \left( _{p+1}F_{p}\right) \left( -t\right)
\right\}  \notag \\
=&1/{\Gamma \left( \alpha +1\right) }\left( _{p+2}F_{p+1}\left( a_{1},\ldots
,a_{p+1},1;b_{1},\ldots ,b_{p},\alpha +1;-t\right) \right).  \label{I-H-p}
\end{align}%
Thus, (\ref{D-H-p}) and (\ref{I-H-p}) show that fractional derivatives and
integrals of the function $\left( _{p+1}F_{p}\right) \left( -t\right) $ are
also analytic in the $t$-plane punctured at the points $0,-1,\infty $. Now
we call attention to the fact that the hypergeometric function $\left(
_{p+1}F_{p}\right) \left( -t\right) $ can be obtained as a result of
successive differentiations and integrations of fractional orders of the
geometric function $\left( 1+t\right) ^{-1}$. Indeed, representing the
geometric function in the form $\left( _{2}F_{1}\left( 1,1;1;-t\right)
\right) $, and using again the relations (\ref{T-ser}) and (\ref{T-coef}),
we have
\begin{equation*}
\left( _{p+1}F_{p}\right) \left( -t\right) \equiv \frac{\Gamma \left(
a_{1}\right) \ldots \Gamma \left( a_{p+1}\right) }{\Gamma \left(
b_{1}\right) \ldots \Gamma \left( b_{p}\right) }\mathcal{D}_{a_{1}-1\ldots }%
\mathcal{D}_{a_{p+1}-1}\mathcal{I}_{b_{1}-1\ldots }\mathcal{I}%
_{b_{p}-1}\left( _{2}F_{1}\left( 1,1;1;-t\right) \right) .
\end{equation*}%
Thus, the function $F\left( t\right) =\left( _{2}F_{1}\left( 1,1;1;-t\right)
\right) $ generates all hypergeometric functions of the form given by (\ref%
{H-p}) and is analytic in the $t$-plane except for the pole at $t=-1$.
Applying formulas (\ref{D-H-p}) and (\ref{I-H-p}) for the function $F\left(
t\right) =\left( _{2}F_{1}\left( 1,1;1;-t\right) \right) ,$ we have%
\begin{equation}
\mathcal{D}_{\alpha }\left\{ F\left( t\right) \right\} =\Gamma \left( \alpha
+1\right) \left( _{2}F_{1}\left( \alpha +1,1;1;-t\right) \right) ,
\label{d-geom}
\end{equation}%
\begin{equation}
\mathcal{I}_{\alpha }\left\{ F\left( t\right) \right\} =1/{\Gamma \left(
\alpha +1\right) }\left( _{2}F_{1}\left( 1,1;\alpha +1;-t\right) \right) .
\label{i-geom}
\end{equation}%
The relation (\ref{d-geom}) can be rewritten in the form
\begin{equation}
F\left( t,\alpha \right) =\mathcal{D}_{\alpha }\left\{ F\left( t\right)
\right\} =\Gamma \left( \alpha +1\right) \left( 1+t\right) ^{-\alpha -1},
\label{be}
\end{equation}%
and if $\alpha $ is not an integer then $\mathcal{D}_{\alpha }\left\{
F\left( t\right) \right\} $ clearly has two singular point at $t=-1$ and $%
t=\infty $. The fractional integral $\mathcal{I}_{\alpha }\left\{ F\left(
t\right) \right\} $ has three singular points at $t=-1,t=\infty ,$ and $t=0.$
In order to verify the last statement we consider the following monodromic
relation for the hypergeometric function in (\ref{i-geom}) where $t^{\ast
}=1-\left( 1+t\right) e^{2\pi i}:$

\begin{equation}
_{2}F_{1}\left( 1,1;\alpha +1;-t^{\ast }\right) -_{2}F_{1}\left( 1,1;\alpha
+1;-t\right) =2\pi ie^{\pi i\alpha }\left( 1+t\right) ^{-1}.
\label{geom-alpha}
\end{equation}%
This relation can be derived from the Euler linear transformation formula.
Letting $t\rightarrow 0$ in (\ref{geom-alpha}) shows that $t=0$ is a
singular point of $\mathcal{I}_{\alpha }\left\{ F\left( t\right) \right\} $.
In fact, such an appearance of new singular points after fractional
differentiation or integration is an obstacle that stopped the creation of
the theory of fractional calculus in the complex plane. However, for the
current case if we continue this process of fractional differentiation or
integration new singularities do not appear. Indeed, applying again (\ref%
{D-H-p}) or (\ref{I-H-p}) to (\ref{d-geom}) we have,%
\begin{equation}
\mathcal{D}_{\beta }\left\{ F\left( t,\alpha \right) \right\} =\Gamma \left(
\beta +1\right) \Gamma \left( \alpha +1\right) \left( _{2}F_{1}\left( \alpha
+1,\beta +1;1;-t\right) \right) ,  \label{D-alpha}
\end{equation}%
\begin{equation}
\mathcal{I}_{\beta }\left\{ F\left( t,\alpha \right) \right\} =\Gamma \left(
\alpha +1\right) /{\Gamma \left( \beta +1\right) }\left( _{2}F_{1}\left(
\alpha +1,1;\beta +1;-t\right) \right) .  \label{I-alpha}
\end{equation}
Functions $\mathcal{D}_{\beta }\left\{ F\left( t,\alpha \right) \right\} $
and $\mathcal{I}_{\beta }\left\{ F\left( t,\alpha \right) \right\} $ are
analytic in the extended $t$-plane punctured at three points $t=-1,t=0$ and $%
t=\infty $. Further differentiation and integration of fractional order do
not change this set of three singular points. Our aim is to explain the
reason for this stabilization and extend this preservation property to wide
classes of \textquotedblleft hypergeometric\textquotedblright\ functions
generated by the system of linear ODE's with analytic coefficients in their
Laplace-Borel dual complex plane.

Previously the topology of the space $S$ was determined by the topology of
the space $H$. In the next section we show how to endow separately both the
spaces $H\left( a,R\right) $ and $S\left( a,R\right) $ by independent
countable normalized topologies.

\subsection{The Basic Spaces and the Duality Theorem}

Given $R\geq 0$ and $a>0,$ we will endow linear spaces $H\left( a,R\right) $
and $S\left( a,R\right) ,$ given by definition 3, by countable normalized
topologies and we retain the previous notation of the normalized spaces.

Thus, $H\left( a,R\right) $ is the set of all functions $F\left( t\right) $
satisfying the conditions:

\begin{itemize}
\item (i) $F\left( t\right) $ is analytic in the region {$\mathfrak{D}$}$(a)$
with boundary $\gamma(a)$ given by (\ref{N-region}) and (\ref{N-boundary}),
respectively;

\item (ii) $F\left( t\right) $ has an exponential growth at infinity;

\item (iii) for every $r>R\geq 0,0<A<a$ the following restriction holds
\begin{equation}
\left\Vert F\right\Vert _{r,A}=\int_{\gamma \left( A\right) }e^{-r\left\vert
t\right\vert }\left\vert F\left( t\right) \right\vert |dt|<\infty .
\label{r-norm}
\end{equation}
\end{itemize}

$S(a,R)$ is a set of all functions $P\left( \zeta \right) ,$ analytic in the
half-plane $\Re \{\zeta \}>R\geq 0$ and satisfying conditions:

\begin{itemize}
\item (i) there exists a sequence $p_{0},p_{1},\ldots ,$ of complex numbers
such%
\begin{equation*}
P\left( \zeta \right) \sim \sum_{k=0}^{\infty }p_{k}/\zeta ^{k},\Re \{\zeta
\}>R,\zeta \rightarrow \infty ;
\end{equation*}

\item (ii) the following sequence of inequalities holds for $0<A<a$ and $%
R<r<\infty $:%
\begin{equation}
\left\vert P\left( \zeta \right) -\sum_{k=0}^{n-1}\frac{p_{k}}{\zeta ^{k}}%
\right\vert \leq \frac{Mn!}{A^{n}\left\vert \zeta \right\vert ^{n}}%
,\,,n=0,1,\ldots ,  \label{P-n-error}
\end{equation}%
where $0<M=M(A,r)<\infty $ is a constant independent of $n;$

\item (iii) for $0<A,R<r<\infty ,\Re \{\zeta \}\geq r$ and $P_{n}\left(
\zeta \right) ,$ defined by (\ref{P-remainder}), the following restriction
holds%
\begin{equation}
\left\Vert P\right\Vert _{r,A}=\sup_{0\leq n<\infty }\frac{A^{n}\left\vert
\zeta \right\vert ^{n}}{n!}\left\vert P_{n}\left( \zeta \right) \right\vert
\leq M(A,r).  \label{P-norm}
\end{equation}
\end{itemize}

We claim the following statement (due to F. Nevanlinna).

\begin{thm}
({Duality Theorem):} $\mathcal{L}H(a,R)=S(a,R)$ and spaces $H(a,R)$ and $%
S(a,R)$ are isomorphic.
\end{thm}

For a further discussion of Nevanlinna's theorems and their extensions, see
\cite{GG}, \cite{Gur}, and references in these papers.

We introduce also the Banach space $H^{1}\left( a\right) $ of functions $%
F\left( t\right) \in H\left( a,R\right) $ that satisfy the following
restriction%
\begin{equation}
\int_{\gamma (a)}\left\vert \frac{F\left( t\right) }{t}\right\vert
|dt|<\infty .  \label{H-1}
\end{equation}

\ \ It follows from (\ref{H-1}) that $F\left( t\right) \in H^{1}\left(
a\right) $ can be represented by the following Cauchy integral%
\begin{equation}
F\left( t\right) =\frac{1}{2\pi i}\int_{\gamma \left( a\right) }\left( 1-%
\frac{t}{\xi }\right) ^{-1}\frac{F\left( \xi \right) }{\xi }d\xi ,
\label{F-Cauchy}
\end{equation}%
which is absolutely convergent for all $t\in ${$\mathfrak{D}$}$(a)$.
Applying operators $\mathcal{D}_{\alpha }$ and $\mathcal{I}_{\alpha }$ to
both the sides of (\ref{F-Cauchy}) and using (\ref{d-geom}) and (\ref{i-geom}%
), we have

\begin{align}
& \mathcal{D}_{\alpha }\left\{ F\left( t\right) \right\} =\frac{\Gamma
\left( \alpha +1\right) }{2\pi i}\int_{\gamma \left( a\right) }\left( 1-%
\frac{t}{\xi }\right) ^{-\alpha -1}\frac{F\left( \xi \right) }{\xi }d\xi ,
\label{alpha-dual-D} \\
& \mathcal{I}_{\alpha }\left\{ F\left( t\right) \right\} =\frac{1}{\Gamma
\left( \alpha +1\right) }\frac{1}{2\pi i}\int_{\gamma \left( a\right)
}\left( _{2}F_{1}\left( 1,1;\alpha +1;\frac{t}{\xi }\right) \right) \frac{%
F\left( \xi \right) }{\xi }d\xi .  \label{alpha-dual-I}
\end{align}%
The functions $\mathcal{D}_{\alpha }\left\{ F\left( t\right) \right\} ,%
\mathcal{I}_{\alpha }\left\{ F\left( t\right) \right\} $ don't belong, in
general, to $H^{1}\left( a\right) .$ However, $\mathcal{D}_{\alpha }\left\{
F\left( t\right) \right\} ,\mathcal{I}_{\alpha }\left\{ F\left( t\right)
\right\} \in H\left( a,0\right) $.

Applying the Laplace transform operator $\mathcal{L}$ to both sides of (\ref%
{alpha-dual-D}) and (\ref{alpha-dual-I}), a straightforward calculations
yield the dual integral representations for
\begin{align}
& P\left( \zeta ,\alpha \right) =\mathcal{\hat{D}}_{\alpha }\left\{ P\left(
\zeta \right) \right\}   \label{P-dual-D} \\
& P^{\ast }\left( \zeta ,\alpha \right) =\mathcal{\hat{I}}_{\alpha }\left\{
P\left( \zeta \right) \right\}   \label{P-dual-I}
\end{align}%
in terms of $F\left( t\right) $, respectively, with kernels--\textit{%
Dingle-Berry basic kernel} and \textit{associated kernel}-- that are given
by (\ref{P-dual-D})and (\ref{P-dual-I}), respectively, for a particular case
$F\left( t\right) =\frac{1}{1+t}$. These kernels will be studied in detail
in the beginning of Section III in Part 2.

\begin{rmk}
The dual integral representation for $P\left( \zeta ,\alpha \right) $ and $%
P^{\ast }\left( \zeta ,\alpha \right) $ in terms of $F\left( t\right) \in
H(a,R),$ exponentially growing at infinity, are very important for
applications to the asymptotic analysis in general, and especially for the
ODE's with analytic coefficients. This approach was initiated in Section 5
of \cite{Gur}, and is further developed in Part 2.
\end{rmk}

In conclusion, we note the important properties of Banach spaces $%
H^{1}\left( a\right) $:

\begin{lmm}
For every $R\geq 0$ the closure of $H^{1}\left( a\right) $ in $H\left(
a,R\right) $ coincides with $H\left( a,R\right) .$
\end{lmm}

In the next sections we will show how to extend the above results given by (%
\ref{alpha-dual-D}) and (\ref{alpha-dual-I}) for general case of
exponentially growing functions of spaces $H(a,R)$.

\subsection{Integral representations for fractional derivatives and integrals%
}

\subsubsection{Integral representations for $\mathcal{\hat{D}}_{\protect%
\alpha }\left\{ P\left( \protect\zeta \right) \right\} $ \ and $\mathcal{%
\hat{I}}_{\protect\alpha }\left\{ P\left( \protect\zeta \right) \right\} $
in terms of $P\left( \protect\zeta \right) $}

\textit{Let }$\Re \{\alpha \}>-1,$ and let $P\left( \zeta \right) \in
S\left( a,R\right) $\textit{\ be given by (\ref{LaplaceF}) for some }$%
F\left( t\right) \in H\left( a,R\right) .$\textit{\ Then using notation (\ref%
{P-calculus}) we can easily derive the validity of the following integral
representation}
\begin{equation}
P\left( \zeta ,\alpha \right) =\mathcal{\hat{D}}_{\alpha }\left\{ P\left(
\zeta \right) \right\} =\frac{\Gamma \left( \alpha +1\right) }{2\pi i}%
\int_{r-i\infty }^{r+i\infty }\left( 1-\frac{z}{\zeta }\right) ^{-\alpha -1}%
\frac{P\left( z\right) }{z}dz,  \label{alpha-D}
\end{equation}%
\textit{where the integral (\ref{alpha-D}) is absolutely convergent for all }%
$\zeta :0\leq R\left( F\right) <\Re \left\{ \zeta \right\} <\infty ,$ and $%
R\left( F\right) $ is given by (\ref{RF}).

We note that this integral representation can be extended to all functions $%
P\left( \zeta \right) $ satisfying (\ref{LaplaceF}). Any additional
constraints on the set of functions $F\left( t\right) $ are not required.
Clearly formula (\ref{alpha-D}) can also be rewritten in the form
\begin{equation}
P\left( \zeta ,\alpha \right) =\frac{\Gamma \left( \alpha +1\right) }{2\pi i}%
\int_{r-i\infty }^{r+i\infty }\left( _{2}F_{1}\left( \alpha +1,1;1;\frac{z}{%
\zeta }\right) \right) \frac{P\left( z\right) }{z}dz.  \label{S-alpha-D}
\end{equation}

A more sophisticated argument shows that for $P\left( \zeta \right) \in
S\left( a,R\right) $ we have%
\begin{equation}
\mathcal{\hat{I}}_{\alpha }\left\{ P\left( \zeta \right) \right\} =\frac{1}{%
2\pi i\Gamma \left( \alpha +1\right) }\int_{r-i\infty }^{r+i\infty }\left(
_{2}F_{1}\left( 1,1;\alpha +1;\frac{z}{\zeta }\right) \right) \frac{P\left(
z\right) }{z}dz,  \label{S-alpha-I}
\end{equation}%
where the integral (\ref{S-alpha-I}) is also absolutely convergent for $\Re
\left\{ \zeta \right\} >R\left( F\right) $. Note that this relation, as
opposed to (\ref{S-alpha-D}), cannot be proved without the additional
restrictions on $F\left( t\right) $ or $P\left( \zeta \right) $. To prove
the validity of (\ref{S-alpha-I}) we must expand both the sides into power
series in $1/\zeta ,$ check that the series are identical and then apply the
uniqueness theorem 1.

We note that the above integral representations can not be used for
applications. Using them we can not even calculate the coefficients of their
asymptotic expansions. In the next section we demonstrate integral
representations for $\mathcal{D}_{\alpha }\left\{ F\left( t\right) \right\} $
and $\mathcal{I}_{\alpha }\left\{ F\left( t\right) \right\}$ in terms of $%
F\left( t\right)$. Their proof is much more complicated, but later we will
show how they can be used in asymptotic analysis of the ODE's to obtain much
simpler final results.

\subsubsection{Integral representations for $\mathcal{D}_{\protect\alpha %
}\left\{ F\left( t\right) \right\} $ and $\mathcal{I}_{\protect\alpha %
}\left\{ F\left( t\right) \right\} $ in terms of $F\left( t\right) $}

Let $H\left( a,R\right) $ and $S\left( a,R\right) $ are our basic spaces
given by definitions 6 and 7, and let operators $\mathcal{D}_{\alpha },%
\mathcal{I}_{\alpha }$ and $\mathcal{\hat{D}}_{\alpha },\mathcal{\hat{I}}%
_{\alpha }$ , are given by (\ref{inH})\ and (\ref{inS}), respectively.

Assume that $0\leq R<\infty ,\,0<A<a<\infty $, and we preserve the notations
$F\left( t,\alpha \right) =\mathcal{D}_{\alpha }\left\{ F\left( t\right)
\right\} $ and $F^{\ast }\left( t,\alpha \right) =\mathcal{I}_{\alpha
}\left\{ F\left( t\right) \right\} ,$ $\Re \{\alpha \}>-1,$\ given by (\ref%
{F-calculus}) and (\ref{P-calculus}), respectively. The following statements
hold:

\begin{thm}
\textit{Let }$F\left( t\right) \in H\left( a,R\right) $\textit{\ and }$%
F(t,\alpha )=\mathcal{D}_{\alpha }\left\{ F\left( t\right) \right\} .$%
\textit{\ Then} \textit{(i) For all }$r,$ $0\leq R<r<\infty $,\textit{\ and}
\textit{for }$t\in $\textit{$\mathfrak{D}$}$(A)$\textit{\ the following
integral representation is valid}
\begin{equation}
F(t,\alpha )=\sum_{k=0}^{\infty }\frac{\Gamma \left( \alpha +k+1\right) }{%
\left( k!\right) ^{2}}\left( rt\right) ^{k}f_{k}(\alpha ,r,t),
\label{F-alpha-R-final}
\end{equation}
\textit{where}
\begin{equation}
f_{k}(\alpha ,r,t)=\frac{1}{2\pi i}\int_{\gamma }\left( _{2}F_{1}\left(
\alpha +k+1,1;k+1;\frac{t}{\xi }\right) \right) \frac{F\left( \xi \right)
e^{-r\xi }}{\xi }d\xi ,  \label{f-k-coef}
\end{equation}
\textit{and the path} $\gamma =\gamma (A)$\textit{\ is oriented in the
counterclockwise direction}.

\textit{(ii) The integral of (\ref{f-k-coef}) and the sum of (\ref%
{F-alpha-R-final}) are absolutely convergent and}
\begin{equation}
F(t,\alpha )\in H\left( a,R\right) .  \label{in B}
\end{equation}
\end{thm}

\begin{thm}
\textit{Let }$F\left( t\right) \in H\left( a,R\right) $\textit{\ and }$%
F^{\ast }(t,\alpha )=\mathcal{I}_{\alpha }\left\{ F\left( t\right) \right\}
. $\textit{\ Then}
\end{thm}

\textit{(i) For all }$r,0\leq R<r<\infty $,\textit{\ and} \textit{for }$t\in
$\textit{$\mathfrak{D}$}$(A),0<A<a,$\textit{\ the following integral
representation is valid}

\ \
\begin{equation}
F^{\ast }(t,\alpha )=\sum_{k=0}^{\infty }\frac{1}{k!\Gamma \left( \alpha
+k+1\right) }\left( rt\right) ^{k}f_{k}^{\ast }(\alpha ,r,t),
\label{F*-alpha-R-final}
\end{equation}%
\textit{where}
\begin{equation}
f_{k}^{\ast }(\alpha ,r,t)=\frac{1}{2\pi i}\int_{\gamma }\left(
_{2}F_{1}\left( k+1,1;\alpha +k+1;\frac{t}{\xi }\right) \right) \frac{%
F\left( \xi \right) e^{-r\xi }d\xi }{\xi },  \label{f*-k-coef}
\end{equation}%
\textit{and the path} $\gamma =\gamma (A)$\textit{\ is oriented in the
counterclockwise direction}.

\textit{(ii) The integral of (\ref{f*-k-coef}) and the sum of (\ref%
{F*-alpha-R-final}) are absolutely convergent and}
\begin{equation}
F^{\ast }(t,\alpha )\in H\left( a,R\right) .  \label{*in B}
\end{equation}

Theorems 2 and 3 together imply Theorem 1.

Proofs of Theorem 6 is given essentially in \cite{Gur}, and Theorem 7 can be
proved using similar technique.

\begin{clr}
\textit{Assume that }$F\left( t\right) \in H\left( a,R\right) $\textit{\
admits an analytical continuation to a region of the Riemann surface of $%
\log t$ that can be represented as a union of overlapping regions of the
form }$D=\cup e^{i\theta ^{\prime }}{\mathfrak{D}}(a^{\prime })$\textit{\
for any set of }$\theta ^{\prime }\in \left( -\infty ,+\infty \right) $%
\textit{\ and }$a^{\prime }\in \left( a,\infty \right). $ \textit{Assume
further that $F\left( t\right)$ retains an exponential growth of the type $R$
in $D$. Then the integral representations given by Theorem 6 and 7 can be
extended to the region }$D$\textit{.}
\end{clr}

\begin{rmk}
\textit{For }$\alpha =0$ \textit{both formulas (\ref{F-alpha-R-final}) and (%
\ref{F*-alpha-R-final}) become trivial identities.}
\end{rmk}

We close this section with the observation that if the function $F\left(
t\right) \in H^{1}\left( a\right) $ then integral representations in
Theorems 6 and 7 simplify considerably. They are given by (\ref{alpha-dual-D}%
) and (\ref{alpha-dual-I}).

\begin{rmk}
Note that the set of all functions $F\left( t\right) \in H^{1}\left(
a\right) $ given by (\ref{H-1}) can be considered as an analog of the Hardy
space $H^{1}\left( D_{a}\right) ,$ where $D_{a}=\left\{ t:\left\vert
t\right\vert <a\right\} .$

Upon completion of this work, we found that similar Hardy spaces in the
regions $\mathfrak{D}(a),$ and in the complement of $\mathfrak{D}(a)$ with
respect to $\mathbb{C}$, were studied by Peschansky (1989), see references
in \cite{Samko}.
\end{rmk}

In what follows, we draw attention to the following result on the analytic
continuation with respect to a parameter $\alpha $\ which has been proved in
\cite{Gur} for the fractional derivatives, and can also be derived for the
fractional integrals.

\subsubsection{Analytic continuation in the $\protect\alpha $-plane}

The following statement is valid.

\begin{thm}
\textit{\ Assume that }$F\left( t\right) $\textit{\ satisfies the
assumptions of Theorem 3. Given fixed }$t\in \mathfrak{D}\left( A\right)
,0<A<a,$\textit{\ }
\end{thm}

\textit{(i) the function }$F(t,\alpha )=\mathcal{D}_{\alpha }\left\{ F\left(
t\right) \right\} $\textit{\ admits an analytic continuation to the whole }$%
\alpha $\textit{-plane except for the negative integers;}

\textit{(ii) the function }$F(t,\alpha )$\textit{\ has a simple pole at
every point }$\alpha =-n-1,n=0,1,\ldots ,$\textit{\ so that the function }
\begin{equation*}
\frac{F(t,\alpha )}{\Gamma \left( \alpha +1\right) }
\end{equation*}%
\textit{is an entire function of exponential type, and the following
relation is valid }
\begin{equation}
\lim_{\alpha \rightarrow -n-1}\frac{F(t,\alpha )}{\Gamma \left( \alpha
+1\right) }=\Psi _{n,F}\left( t\right) ,  \label{B-polyn}
\end{equation}%
\textit{\ where }$\Psi _{n,F}\left( t\right) $\textit{\ is a polynomial in }$%
t$\textit{\ of degree }$n,$\textit{\ which can be represented in the form}
\begin{equation}
\Psi _{n,F}\left( t\right) =\sum_{j=0}^{n}\left( -1\right) ^{j}\left(
\begin{array}{c}
n \\
j%
\end{array}%
\right) f_{j}t^{j},  \label{Bern-r}
\end{equation}%
\textit{with coefficients} $f_{j}={F^{(j)}}(0)/j!,$\textit{\ which can be
represented as}
\begin{equation}
f_{j}=\frac{1}{2\pi i}\sum_{s=0}^{j}\frac{1}{\left( j-s\right) !}%
r^{j-s}\int_{\gamma \left( A\right) }\frac{F\left( \xi \right) e^{-r\xi
}d\xi }{\xi ^{1+s}},  \label{coefficients}
\end{equation}

\textit{(iii) }$F(t,\alpha )=\mathcal{I}_{\alpha }\left\{ F\left( t\right)
\right\} $\textit{\ is an entire function of an exponential type.}

We note that Theorems 6--8 extend the corresponding results in \cite{Osler1}
and \cite{Osler2}, that was described also in \cite{Samko}, to our version
of fractional calculus.

In the next section we show how to apply the above version of fractional
calculus to linear ODE's with analytic coefficients.

\section{Linear ODE's with analytic coefficients}

In this section we study linear spaces of vector-valued functions generated
by the second order linear ODE's with analytic coefficients in their
Laplace-Borel dual complex plane. We discovered that the elements of these
spaces inherit a number of properties of classical hypergeometric functions.
In particular, we show how to extend some classical linear transformation
formulae for hypergeometric functions (due to Euler, Gauss, Goursat) to the
elements of the dual spaces. For the cases of the Kummer or Whittaker
differential equations the above elements can be represented in terms of
classical hypergeometric function. Even for this special case our formulas
are new. In Part II we show how to apply these formulas in asymptotic
analysis of the ODE's.

\subsection{Introduction}

We develop further an approach to study asymptotic properties of solutions
of linear ODE's with coefficients that are analytic in a neighborhood of
infinity first proposed in \cite{GG-1} for the case of even coefficients.
The approach is based on the duality between a linear system of monodromic
functional equations generated by the ODE and its Laplace-Borel-dual system.
Here we consider the case of a second order ODE, which generates a system of
two monodromic integral equations of a fractional order for a pair of
multi-valued analytic functions: a jump of one function across the cut is
expressed in terms of the fractional integral of the other function.
Moreover, we stated for this pair of functions a system of Euler and
Euler-Goursat type linear transformation formulae. Our approach is based on
a version of fractional calculus in the complex plane that was initiated in
\cite{Gur}. It can be extended to wide classes of matrix-valued function.

\subsection{A system of functional monodromic equations generated by the ODE
in the Laplace-Borel dual complex plane}

We consider linear ODE's of the form
\begin{equation}
u^{\prime \prime }(\zeta )+a(\zeta )u^{\prime }(\zeta )+b(\zeta )u(\zeta )=0,
\label{Kummer-ODE}
\end{equation}%
where functions $a(\zeta )$ and $b(\zeta )$ are analytic at infinity. Set $%
a(\infty )=a_{0}$ and $b\left( \infty \right) =b_{0}$. If $a_{0}-4b_{0}\neq
0 $ then this equation can be reduced to a perturbation of the Whittaker
differential equation (pWde):
\begin{equation}
\frac{d^{2}u}{d\zeta ^{2}}=\left( 1/4-\kappa /\zeta +\left( \mu
^{2}-1/4\right) /\zeta ^{2}+\left( 1/\zeta ^{3}\right) \sum_{k=0}^{\infty
}\left( \beta _{k}/\zeta ^{k}\right) \right) u,  \label{p-even-WH}
\end{equation}%
where $\kappa ,\mu ,\beta _{k},k=0,1,\ldots $ are complex numbers and the
series $\sum_{k=0}^{\infty }\left( \beta _{k}/\zeta ^{k}\right) $ is
absolutely convergent in the exterior of a circle of radius $R$ centered at
the origin for some $R\geq 0$.

There exists a pair of linearly independent solutions of (\ref{p-even-WH}), $%
u_{1}\left( \zeta \right) $ and $u_{2}\left( \zeta \right) $, such that
\begin{align}
& u_{1}\left( \zeta \right) =e^{-\frac{1}{2}\zeta }\zeta ^{\kappa
}P_{1}\left( \zeta \right) ,  \label{1} \\
& u_{2}\left( \zeta \right) =e^{\frac{1}{2}\zeta }\zeta ^{-\kappa
}P_{2}\left( \zeta \right) ,  \label{2}
\end{align}%
where \textit{the phase-amplitudes} $P_{1}\left( \zeta \right) $ and $%
P_{2}\left( \zeta \right) $ satisfy the following relations
\begin{align}
& P_{1}\left( \zeta \right) =1+o\left( 1\right) ,\,-\pi \leq \arg \zeta \leq
\pi ,\left\vert \zeta \right\vert >R,  \label{P1} \\
& P_{2}\left( \zeta \right) =1+o\left( 1\right) ,\,0\leq \arg \zeta \leq
2\pi ,\left\vert \zeta \right\vert >R.  \label{P2}
\end{align}%
We note that the solutions $u_{1}\left( \zeta \right) $ and $u_{2}\left(
\zeta \right) $ are uniquely determined by the conditions (\ref{P1}) and (%
\ref{P2}), respectively.

Clearly, $P_{1}\left( \zeta \right) $ and $P_{2}\left( \zeta \right) $ are
analytic multi-valued functions in the exterior of the circle $|\zeta |=R$
in the $\zeta $-plane, and it can be proved that there exists a pair of
complex numbers $T_{1}, T_{2}$ such that the following system of monodromic
functional relations is valid for $\zeta >R:$
\begin{eqnarray}
P_{1}\left( \zeta e^{\pi i}\right) -P_{1}\left( \zeta e^{-\pi i}\right)
&=&T_{1}e^{-\zeta }\zeta ^{-2\kappa }P_{2}\left( \zeta e^{\pi i}\right) ,
\label{MW1} \\
P_{2}\left( \zeta e^{\pi i}\right) -P_{2}\left( \zeta e^{-\pi i}\right)
&=&T_{2}e^{\zeta }\zeta ^{2\kappa }P_{1}\left( \zeta e^{-\pi i}\right) ,
\label{MW2}
\end{eqnarray}%
Thus, the pWde, with an infinite number of parameters, gives rise to a
system of monodromic relations with only three parameters $T_{1}, T_{2}$ and
$\kappa.$

The terminology \textit{monodromic} means that a jump of one multi-valued
function across the cut along any ray $l_{\theta }=\left\{ \zeta :\arg \zeta
=\theta ,\left\vert \zeta \right\vert >R\right\} ,-\infty <\theta <\infty ,$
can be expressed in terms of the other function. The constants
\begin{equation}
T_{1}=T_{1}\left( \kappa ,\mu ,\beta _{0},\ldots \right) ,T_{2}=\left(
\kappa ,\mu ,\beta _{0},\ldots \right)  \label{T-T}
\end{equation}%
are invariants of (\ref{p-even-WH}) usually referred to as the \textit{%
connection coefficients}, or \textit{Stokes multipliers}. The parameter $%
\kappa ,$ which is the second coefficient in (\ref{p-even-WH}), can be
expressed in terms of $T_{1}$ and $T_{2}$.

\begin{dfn}
Using the complex numbers $T_{1}$, $T_{2}$ and $\kappa ,$ we introduce the
linear space $S=S\left( T_{1},T_{2},\kappa \right) $ of all vector-valued
solutions $\left( P_{1}\left( \zeta \right) ,P_{2}\left( \zeta \right)
\right) $ of the system of functional equations (\ref{MW1})--(\ref{MW2}),
which retain the analytic properties of $P_{1}\left( \zeta \right) $ and $%
P_{2}\left( \zeta \right) $ described above.
\end{dfn}

\subsection{Dual system of monodromic relations}

Assume that the vector-valued function $\left( P_{1}\left( \zeta \right)
,P_{2}\left( \zeta \right) \right) $ belongs to the linear space $S=S\left(
T_{1},T_{2},\kappa \right) $ given by Definition . Then there exists a pair
of integrable functions $F_{1}\left( t\right) $ and $F_{2}\left( t\right) $
such
\begin{equation}
P_{1}\left( \zeta \right) =\mathcal{L}\left\{ F_{1}\left( t\right) \right\}
,P_{2}\left( \zeta \right) =\mathcal{L}\left\{ F_{2}\left( t\right) \right\}
,  \label{dual functions}
\end{equation}%
where $\mathcal{L}$ is the Laplace transform operator given by (\ref%
{LaplaceF}). Indeed,it follows from (\ref{P1}) that $P_{1}\left( \zeta
\right) /\zeta $ is a square integrable function on every line $\Im \left\{
\zeta \right\} >R\geq 0$. It follows from Parseval's theorem that $%
F_{1}\left( t\right) $ is a square integrable function on the interval $%
\left( R,+\infty \right) .$ Similar argument can be applied to $P_{2}\left(
\zeta \right) $.

\begin{dfn}
We denote by $H\left( T_{1},T_{2},\kappa \right) $\ the linear space of all
vector-valued functions $\left( F_{1}\left( t\right) ,F_{2}\left( t\right)
\right) $ given by (\ref{dual functions}), so that%
\begin{equation}
H\left( T_{1},T_{2},\kappa \right) =\mathcal{L}S\left( T_{1},T_{2},\kappa
\right) ,  \label{H=LS}
\end{equation}
\end{dfn}

In what follows we study properties of elements of $H\left(
T_{1},T_{2},\kappa \right) $. Applying the Borel transform operator $%
\mathcal{L}^{-1}$ to the left and right sides of (\ref{MW1})--(\ref{MW2}),
we show that the jumps of $P_{1}\left( \zeta \right) $ and $P_{2}\left(
\zeta \right) $ on the left sides are transformed into the jump of $%
F_{1}\left( t\right) $ and $F_{2}\left( t\right) $ on the corresponding cut,
while the exponential factors and power factors on the right sides are
transformed into the shifts and the fractional integrals of orders of $%
-2\kappa $ and $2\kappa $ of $F_{2}\left( t\right) $ and $F_{1}\left(
t\right) $, respectively. The latter fact follows from Theorem 4. This
enables us to derive the dual system of monodromic relations for a pair $%
F_{1}\left( t\right) $ and $F_{2}\left( t\right) .$

Now we demonstrate properties of $F_{1}\left( t\right) $ and $F_{2}\left(
t\right) $ that follows from (\ref{P1}) and (\ref{P2}) and (\ref{MW1})--(\ref%
{MW2}).

Assume that $F_{1}\left( t\right) $\ and $F_{2}\left( t\right) $ are given
by (\ref{dual functions})\bigskip . Then the following statements are valid

\begin{thm}
\begin{itemize}
\item (i) $F_{1}\left( t\right) $\ and $F_{2}\left( t\right) $ are\ analytic
in the $t$-plane cut along the intervals $\left( -\infty ,-1\right) $ and $%
\left( 1,+\infty \right) $, respectively;

\item (ii) $F_{1}\left( t\right) $\ and $F_{2}\left( t\right) $ can be
continued analytically further to the Riemann surface of $\log t,$ oriented
in the counterclockwise direction, except for the points $\left( -1,0\right)
,\left( 0,0\right) $ and $\left( 0,0\right) ,\left( 0,1\right) $,
respectively;

\item (iii) In every sectorial region of the form%
\begin{equation*}
S\left( \theta _{1},\theta _{2}\right) =\left\{ t:\theta _{1}<\arg t<\theta
_{2},0<\left\vert t\right\vert <\infty \right\}
\end{equation*}
$F_{1}\left( t\right) $\ and $F_{2}\left( t\right) $ have an exponential
growth of type less or equal to $R$ at infinity.
\end{itemize}
\end{thm}

\begin{thm}
The following system of monodromic relations%
\begin{eqnarray}
F_{1}\left( te^{-\pi i}\right) -F_{1}\left( te^{\pi i}\right) &=&T_{1}\left(
t-1\right) ^{2\kappa }\mathcal{I}_{2\kappa }\left\{ F_{2}\left( \left(
t-1\right) e^{-\pi i}\right) \right\} ,  \label{Mon-1} \\
F_{2}\left( te^{-\pi i}\right) -F_{2}\left( te^{\pi i}\right) &=&T_{2}\left(
te^{\pi i}-1\right) ^{-2\kappa }\mathcal{I}_{-2\kappa }\left\{ F_{2}\left(
\left( te^{\pi i}-1\right) \right) \right\} ,  \label{Mon-2}
\end{eqnarray}%
is valid in the whole $t$-plane punctured in the points $(-1,0),(0,0)$ and $%
(0,1).$
\end{thm}

The system of relations (\ref{Mon-1})--(\ref{Mon-2}), which is dual to the
system (\ref{MW1})--(\ref{MW2}), can be derived from the latter system
directly applying to both sides of it the operator $\mathcal{L}^{-1}$ ,
using the relation (\ref{alt-I}) with $\alpha =\pm {2\kappa }.$

\begin{dfn}
Using the complex numbers $T_{1}$, $T_{2}$ and $\kappa $ of the Definition
3, we introduce the linear space $\bar{H}\left( T_{1},T_{2},\kappa \right) $
of all vector-valued solutions $\left( F_{1}\left( t\right) ,F_{2}\left(
t\right) \right) $ of the system of functional equations (\ref{Mon-1})--(\ref%
{Mon-2}), where $F_{1}\left( t\right) $\ and $F_{2}\left( t\right) $ satisfy
the properties (i)--(iii) of Theorem 9.
\end{dfn}

\begin{rmk}
We note that the linear space $\bar{H}\left( T_{1},T_{2},\kappa \right) $\
is essentially wider then the linear space $H\left( T_{1},T_{2},\kappa
\right) $ given by (\ref{H=LS} )$.$ This fact was noted in paper \cite{GG-1}
where we considered the case when the coefficients of (\ref{Kummer-ODE}) are
even function. For this case $T_{1}=T_{2}=T,\kappa=0,$ and the linear space $%
S\left( T_{1},T_{2},\kappa \right) $ was denoted by $S_{1,T}$. See then
Definition 1 and the last line of the second paragraph after on page 658 of
\cite{GG-1}. Therefore it is impossible in general to derive the system (\ref%
{MW1})--(\ref{MW2}) applying the Laplace transform operator to the system (%
\ref{Mon-1})--(\ref{Mon-2}) without imposing addition conditions on the
functions $F_{1}\left( t\right) $\ and $F_{2}\left( t\right). $
\end{rmk}

The question arises whether there is a system of functional equations for
functions $F_{1}\left( t\right) $\ and $F_{2}\left( t\right), $ satisfying
the properties (i)--(iii), which is equivalent to the system (\ref{MW1})--(%
\ref{MW2}) and which entails the system (\ref{Mon-1})--(\ref{Mon-2})?

Now we see that the answer to this question is affirmative.

\subsection{Euler-Gauss-type system of linear transformation formulas}

\begin{thm}
Let $F_{1}\left( t\right) $ and $F_{2}\left( t\right) $ are given by (\ref%
{dual functions}) where $\left( P_{1}\left( \zeta \right) ,P_{2}\left( \zeta
\right) \right) \in $. Then the pair $\left( F_{1}\left( t\right)
,F_{2}\left( t\right) \right) $, if $2\kappa $ is not integer, satisfies a
system of Euler-Gauss-type linear transformation formulae of the form%
\begin{align}
& 2i\sin \left( 2\kappa \pi \right) F_{1}\left( t\right) =-T_{1}\left(
1+t\right) ^{2\kappa }\mathcal{I}_{2\kappa }\left\{ F_{2}\left( 1+t\right)
\right\}  \notag \\
& +T_{2}e^{-2\kappa \pi i}\mathcal{I}_{-2\kappa }\left\{ F_{1}\left( \left(
1+t\right) e^{-\pi i}\right) \right\} ,-\pi \leq \arg t\leq \pi
,1<|t|<\infty ,  \label{EG-LTF-1}
\end{align}%
\begin{align}
& 2i\sin \left( 2\kappa \pi \right) F_{2}\left( t\right) =T_{2}\left(
t-1\right) ^{-2\kappa }\mathcal{I}_{-2\kappa }\left\{ F_{1}\left( t-1\right)
\right\}  \notag \\
& -T_{1}\mathcal{I}_{2\kappa }\left\{ F_{2}\left( \left( t-1\right) e^{-\pi
i}\right) \right\} ,-2\pi \leq \arg t\leq 0,1<|t|<\infty ,  \label{EG-LTF-2}
\end{align}%
where $\mathcal{I}_{\alpha }$ is the operator of a fractional integration of
order $\alpha =\pm 2\kappa $ in our version of fractional calculus in the
complex plane, given by (\ref{inH}).
\end{thm}

We recall that, as follows from Theorem 8 (iii), the function $\mathcal{I}%
_{2\kappa }\left\{ F_{j}\left( t\right) \right\} ,j=1,2,$ for fixed $t$, is
entire function of $\kappa .$

\subsection{Euler-Goursat-type linear transformation formulas}

\begin{thm}
If $2\kappa \in \mathbb{Z}$ then the system of an Euler-Goursat-type linear
transformation formulae can be written in the form%
\begin{align}
& F_{1}\left( t\right) =C_{1}T_{1}\left( 1+t\right) ^{2\kappa }\log \left(
1+t\right) \mathcal{I}_{2\kappa }\left\{ F_{2}\left( 1+t\right) \right\}
\notag \\
& +\Psi _{1}\left( 1+t\right) ,-\pi \leq \arg t\leq \pi ,1<|t|<\infty ,
\label{E-Go-1}
\end{align}%
\begin{align}
& F_{2}\left( t\right) =C_{2}T_{2}\left( t-1\right) ^{-2\kappa }\log \left(
t-1\right) \mathcal{I}_{-2\kappa }\left\{ F_{1}\left( t-1\right) \right\}
\notag \\
& +\Psi _{2}\left( t-1\right) ,-2\pi \leq \arg t\leq 0,1<|t|<\infty ,
\label{E-Go-2}
\end{align}%
where $C_{1}$ and $C_{2}$ are constants, which can be calculated exactly,
functions $\Psi _{1}\left( t\right) $ and $\Psi _{2}\left( t\right) $ are
analytic in the unit circle of the $t$-plane, and their Taylor coefficients
can be found by recurrence following the procedure, which was described in
\cite{GG-1} for the case $\kappa =0,\beta _{2k}=0,k=0,1,\ldots $.
\end{thm}

\subsection{Duality theorem}

Our main result, the duality theorem:

\begin{thm}
For the case when $2\kappa \notin \mathbb{Z}$ the linear space $H\left(
T_{1},T_{2},\kappa \right) $ given by (\ref{H=LS}) coincides with the set of
all vector-valued solutions $\left( F_{1}\left( t\right) ,F_{2}\left(
t\right) \right) $ of the system of functional equations (\ref{EG-LTF-1})--(%
\ref{EG-LTF-2}), which satisfy conditions (i) and (ii). For the case when $%
2\kappa \in \mathbb{Z}$ the linear space $H\left( T_{1},T_{2},\kappa \right)
$ coincides with the set of all vector-valued solutions $\left( F_{1}\left(
t\right) ,F_{2}\left( t\right) \right) $ of the system of functional
equations (\ref{E-Go-1})--(\ref{E-Go-2}), which satisfy conditions
(i)--(iii) of Theorem 9.
\end{thm}

\subsection{The case of the Whittaker equation}

Validity of all above result result can be verified using the case of the
(unperturbed) Whittaker equation for which $F_{1}\left( t\right) $\ and $%
F_{2}\left( t\right) $ are the hypergeometric functions%
\begin{align}
& F_{1}\left( t\right) =\left( _{2}F_{1}\left( \frac{1}{2}-\kappa -\mu ,%
\frac{1}{2}-\kappa +\mu ;1;-t\right) \right) ,  \label{M-W1} \\
F_{2}\left( t\right) & =\left( _{2}F_{1}\left( \frac{1}{2}+\kappa -\mu ,%
\frac{1}{2}+\kappa +\mu ;1;-\left( te^{\pi i}\right) \right) \right) ,
\label{M-W2}
\end{align}%
and for the case when $2\kappa \notin \mathbb{Z}$ the relations (\ref%
{EG-LTF-1})--(\ref{EG-LTF-2}) can be derived from Euler's linear
transformation formulae (see \cite{Abramowitz}, \textbf{15.3.6}) written
separately for two functions $F_{1}\left( t\right) $ and $F_{2}\left(
t\right) $ given by (\ref{dual functions}) and rewritten as (\ref{M-W1}) and
(\ref{M-W2}). For the case when $2\kappa \in \mathbb{Z}$ the relations (\ref%
{E-Go-1})--(\ref{E-Go-2}) can be derived from the Euler-Goursat linear
transformations formulae (see \cite{Abramowitz}, \textbf{15.3.10}--\textbf{%
15.3.14}).

\subsection{Euler's linear transformation formula for $\left(
_{2}F_{1}\left( a,b;c;t\right) \right) $}

The hypergeometric function $_{2}F_{1}\left( a,b;c;t\right) $ can be
introduced as an analytic continuation of the hypergeometric series%
\begin{equation}
\frac{\Gamma \left( c\right) }{\Gamma \left( a\right) \Gamma \left( b\right)
}\sum_{k=0}^{\infty }\frac{\Gamma \left( a+k\right) \Gamma \left( b+k\right)
}{\Gamma \left( c+k\right) k!}t^{k},\left\vert t\right\vert <1,
\label{H-series}
\end{equation}%
from the unit circle, along any path not crossing the ponts $t=0$ and $t=1$,
to the $t$-plane.

Setting
\begin{equation*}
A:=\frac{\Gamma \left( c\right) \Gamma \left( a+b-c\right) }{\Gamma \left(
a\right) \Gamma \left( b\right) }\text{ and }B:=\frac{\Gamma \left( c\right)
\Gamma \left( c-a-b\right) }{\Gamma \left( c-a\right) \Gamma \left(
c-b\right) },
\end{equation*}%
{Euler's linear transformation formula,}
\begin{align}
& _{2}F_{1}\left( a,b;c;t\right) =\left( 1-t\right) ^{c-a-b}A\left(
_{2}F_{1}\left( c-a,c-b;c-a-b+1;1-t\right) \right)  \label{Euler's-LTF} \\
& +B\left( _{2}F_{1}\left( a,b;a+b-c+1;1-t\right) \right) ,  \notag
\end{align}%
for hypergeometric functions is valid if
\begin{equation}
c-a-b\notin \mathbb{Z}\text{ and }\left\vert \arg \left( 1-t\right)
\right\vert \leq \pi .  \label{LTF-conditions}
\end{equation}%
This formula is available in all manuals on special functions and still is
in the focus of many mathematicians, see recent publications, \cite{Beukers}%
, \cite{Heckman}, \cite{Kuznetsov}.

Analysis of (\ref{Euler's-LTF}) shows that $_{2}F_{1}\left( a,b;c;t\right) $
given by (\ref{H-series}) admits an analytic continuation to the $t$-plane
cut along the interval $\left( 1,\infty \right) $ of the positive ray and
restriction $\left\vert \arg \left( 1-t\right) \right\vert =\pm \pi $ means
that $t$ belongs to the boundary of the cut $t$-plane. In the next section
we calculate the jump of the hypergeometric function $\left( _{2}F_{1}\left(
a,b;c;t\right) \right) $ on the cut $\left( 1,+\infty \right) .$ But first,
we introduce two transformations
\begin{align}
& \mathcal{K}_{_{1}}\left\{ _{2}F_{1}\left( a,b;c;t\right) \right\} =\left(
_{2}F_{1}\left( c-a,c-b;c-a-b+1;1-t\right) \right) ,  \label{K-1-2} \\
& \mathcal{K}_{_{2}}\left\{ _{2}F_{1}\left( a,b;c;t\right) \right\} =\left(
_{2}F_{1}\left( a,b;a+b-c+1;1-t\right) \right) ,  \notag
\end{align}%
The relations (\ref{Euler's-LTF}) does not provide any clue of how operators
$\mathcal{K}_{_{1}}$\ and $\mathcal{K}_{_{2}}$ are connected. We answer this
question not only for hypergeometric functions given by (\ref{M-W1}) and (%
\ref{M-W2}) but also for wide classes of \textquotedblleft
hypergeometric\textquotedblright\ functions that are generated by solutions
of the system (\ref{MW1})--(\ref{MW2}) in the Laplace-Borel dual complex
plane. It turns out that the above transformations are operators of
fractional integration of order $\kappa $ and $-\kappa $\ in the author's
version of the fractional calculus in the complex plane. In the next Section
we explain how to derive the system of linear transformation formulas for
the functions (\ref{MW1})--(\ref{MW2}) of the form given by (\ref{EG-LTF-1}%
)--(\ref{EG-LTF-2}) from the formulas (\ref{Euler's-LTF}) written separately
for (\ref{MW1})--(\ref{MW2}).

\subsection{Monodromic relations and connection coefficients for $\left(
_{2}F_{1}\left( a,b;c;t\right) \right) $}

Assume that $t$ is any point inside the circle $\left\vert 1-t\right\vert
<1, $ which does not belong to its radius along the interval $\left(
1,2\right) $ of the real line. Since $1-t^{\ast }=\left( 1-t\right) e^{\pm
2\pi i},$ the point $t^{\ast }$ satisfies the same restriction. Let us write
Euler's linear transformation formulas (\ref{Euler's-LTF}) for $%
_{2}F_{1}\left( a,b;c;t\right) $ and $_{2}F_{1}\left( a,b;c;t^{\ast }\right)
$ and note that both the hypergeometric functions in the right hand-side of (%
\ref{Euler's-LTF}) are coincide at $t^{\ast }$and $t$. Then subtracting the
first relation from the second relation, a straightforward calculation
yields the following result
\begin{align}
& _{2}F_{1}\left( a,b;c;t^{\ast }\right) -\left( _{2}F_{1}\left(
a,b;c;t\right) \right)  \notag \\
& =T^{\pm }(a,b,c)\left( 1-t\right) ^{c-a-b}\left( _{2}F_{1}\left(
c-a,c-b;c-a-b+1;1-t\right) \right) ,  \label{GM}
\end{align}
where
\begin{equation}
T^{\pm }(a,b,c)=\mp 2{\pi }ie^{\pm \pi i\left( c-a-b\right) }\frac{\Gamma
\left( c\right) }{\Gamma \left( a\right) \Gamma \left( b\right) \Gamma
\left( c-a-b+1\right) }.  \label{Connec}
\end{equation}
As oppose to the relation (\ref{Euler's-LTF}) both relations (\ref{GM}) and (%
\ref{Connec}) are valid for the case $a+b-c\in \mathbb{Z},$ and can be
extended for all $t$ belonging to the closure of the $t$-plane cut along the
interval $\left( 1,+\infty \right) $. We refer \ to the relation (\ref{GM})
and the constant $T^{\pm }(a,b,c)$ as the \textit{monodromic relation} and
\textit{connection coefficient, }respectively\textit{.}

In the cases when $a+b-c\in \mathbb{Z}$ the formulas (\ref{Euler's-LTF}) are
meaningless while the formulas (\ref{GM}) and (\ref{Connec}) remain valid.
However, the latter formulas do not allow to recognize the nature of the
singularity of the hypergeometric function at $t=1.$ Goursat discovered that
the right hand side of (\ref{Euler's-LTF}) should be replaced by a sum of
two terms, where the first term is essentially the same as the first term in
(\ref{Euler's-LTF}) but the factor $\left( 1-t\right) ^{c-a-b}$ is replaced
by $\left( 1-t\right) ^{c-a-b}\cdot \log \left( 1-t\right) ,$ while the
second term is a function that are analytic in the circle $\left\vert
1-t\right\vert <1$ and its Taylor coefficients were provided with explicit
expression. The corresponding formulas are given in \cite{Abramowitz},
\textbf{15.3.10--15.3.14}. Following ideas of \cite{GG-1}, we generalized
these formulas for the classes of \textquotedblleft
hypergeometric\textquotedblright\ functions, generated by solutions of the
system of monodromic relations in the $t$-plane that is dual to the system (%
\ref{MW1})--(\ref{MW2}).

We note that the relations (\ref{Euler's-LTF}), (\ref{GM}) and (\ref{Connec}%
) can be rewritten as
\begin{align}
& _{2}F_{1}\left( a,b;c;-t\right) =\left( 1+t\right) ^{c-a-b}A\left(
_{2}F_{1}\left( c-a,c-b;c-a-b+1;1+t\right) \right)  \label{LTF-t} \\
& +B\left( _{2}F_{1}\left( a,b;a+b-c+1;1+t\right) \right) ,  \notag
\end{align}%
\begin{align}
& _{2}F_{1}\left( a,b;c;-t^{\ast }\right) -\left( _{2}F_{1}\left(
a,b;c;-t\right) \right)  \notag \\
& =T^{\pm }(a,b,c)\left( 1+t\right) ^{c-a-b}\left( _{2}F_{1}\left(
c-a,c-b;c-a-b+1;1+t\right) \right) ,  \label{GM-t}
\end{align}%
where $1+t^{\ast }=\left( 1+t\right) e^{\pm 2\pi i}$ and \
\begin{equation}
T^{\pm }(a,b,c)=\mp 2{\pi }ie^{\pm \pi i\left( c-a-b\right) }\frac{\Gamma
\left( c\right) }{\Gamma \left( a\right) \Gamma \left( b\right) \Gamma
\left( c-a-b+1\right) }.  \label{Connec-t}
\end{equation}

\subsection{Conclusion}

In conclusion we note that the above systems of linear transformation
formulae and monodromic relations can be extended to the more general case
of a matrix differential equation of the form $du/d\zeta =Q\left( \zeta
\right) u,$ where $u=u\left( \zeta \right) $ is an $n\times n$ matrix-valued
function, $Q\left( \zeta \right) =A_{0}+\frac{1}{\zeta }A_{1}+\ldots ,$
where $A_{0},A_{1},\ldots $ are the constant $n\times n$ matrices, and $\det
A_{0}\neq 0,A_{0}=\diag{\left( \lambda _{1},\ldots \lambda _{n}\right)}%
,\lambda _{i}\neq \lambda _{j},1\leq i,j\leq n,$ $A_{1}=%
\diag{\left( \kappa
_{1},\ldots \kappa _{n}\right)},$ and $Q\left( \zeta \right) $ is analytic
in a neighborhood of infinity. Then associated system of monodromic
relations is expressed in terms of fractional integrals of orders $\pm
\left( \kappa _{1}-\kappa _{2},\ldots ,\kappa _{n-1}-\kappa _{n},\kappa
_{n}-\kappa _{1}\right) .$

{Department of Mathematics, Swinburne University of Technology, PO Box 218,
Hawthorn, Victoria 3122, Australia,\newline
e-mail: vgurarii@swin.edu.au.}

\end{document}